\documentclass[a4paper,12pt]{article}
\usepackage{amsmath, amsfonts, amssymb}
\usepackage[english, russian]{babel}
\begin{document}
\textit{УДК 517.956.6; 517.44}
\begin{center}
\textbf{Задача Коши для вырождающегося гиперболического уравнения
второго рода}\\
\medskip
\textbf{\verb"Т.Г.Эргашев"}
\smallskip \textit{ertuhtasin@mail.ru}
\end{center}

           Maqolada    spektral    parametrli   ikkinchi    tur    buziladigan
           giperbolik   tenglama   uchun   Koshi   masalasi   Riman   usuli
           yordamida  o'rganilgan.    Bu   yerda   kiritilgan   ikki  va   uch
           o'zgaruvchili yangi konflyuent gipergeometrik funksiyalarning
           ba'zi  xossalarini  isbotlash  yo'li  bilan  qo'yilgan  masalaning
           yagona yechimi oshkor ko'rinishda topilgan.

           In this article the  unique  solution  of  the  Cauchy  problem  is
           founded by the Riemann method. Some relations for given here
           confluent hypergeometric functions of  two  and  three  variables
           are used.

     Рассмотрим уравнение
                                    $$x^nV_{xx}-(-y)^mV_{yy}+\mu V=0\eqno(1)$$
в конечной односвязной области $D$, ограниченной характеристиками
$A_0C:\,\,\,\xi =0$, $A_1C:\,\,\,\eta =1$ и $A_0A_1:\,\,\,\,y=0$
уравнения (1) при $y\le 0$ и $x\ge 0$, где $m,$ $n$ и $\mu $
действительные числа, причем $0<m<1,$ $0\le n<1;$
                                 $$\left. \begin{matrix}
   \xi   \\
   \eta   \\
\end{matrix} \right\}=\frac{2}{2-n}{x^{\frac{2-n}{2}}}\pm \frac{2}{2-m}{{(-y)}^{\frac{2-m}{2}}}\eqno(2)$$

З а д а ч а  К о ш и  для уравнения (1) состоит в определении
функции $V(x,y)\in {{C}^{2}}(D)\cap C(\bar{D})$, удовлетворяющей в
области $D$ уравнению (1) и начальным условиям
$$V(x,0)=\tau _1(x),\,\,\, V_y(x,0)=\nu_1(x),\eqno(3)$$
где $\tau_1$,$\nu_1$ - заданные функции, причем $\tau_1\in
C^3[0,1]$,$\nu_1\in C^2[0,1]$.

В характеристических координатах $\xi $ и $\eta $ уравнение (1)
переходит в уравнение типа обобщенного уравнения
Эйлера-Пуассона-Дарбу
 $$u_{\xi \eta }+\frac{\alpha }{\eta +\xi }\left( {{u}_{\eta }}+u_\xi \right)-\frac{\beta }{\eta -\xi }\left( {{u}_{\eta }}-{{u}_{\xi }} \right)-{{\lambda }^{2}}u=0,\eqno(4)$$
а область $D$ преобразуется в область $\Delta $, граница которой
состоит из отрезков прямых $PM:\xi =0$, $QM:\eta =1$ и $PQ:\eta
=\xi $, а начальные условия принимают вид $$u(\xi ,\xi )=\tau (\xi
), 0\le \xi \le 1,\eqno(5) $$ $${{\left[ 2(1-2\beta )
\right]}^{-2\beta }}\underset{\eta \to \xi }{\mathop{\lim
}}\,{{(\eta -\xi )}^{2\beta }}\left( {{u}_{\xi }}-{{u}_{\eta }}
\right)=\nu (\xi ), 0<\xi <1, \eqno(6)$$ где $$2\alpha
=\frac{n}{n-2}, 2\beta =\frac{m}{m-2}, -1<2\alpha \le 0, -1<2\beta
<0, {{\lambda }^{2}}=-\frac{1}{4}\mu ,$$ $$\tau (\xi )={{\tau
}_{1}}\left[ {{\left( \frac{2-n}{2}\xi \right)}^{\frac{2}{2-n}}}
\right], \nu (\xi )={{\nu }_{1}}\left[ {{\left( \frac{2-n}{2}\xi
\right)}^{\frac{2}{2-n}}} \right],$$ $$u(\xi ,\eta )=V\left[
{{\left( \frac{2-n}{4}(\eta +\xi )
\right)}^{\frac{2}{2-n}}},\,\,-\,\,{{\left( \frac{2-m}{4}(\eta
-\xi ) \right)}^{\frac{2}{2-m}}} \right].$$

Для уравнения (1) в случае вырождения первого рода изучены задача
Коши [1] и задачи со смещением [2,3]. Для вырождающихся уравнений
гиперболического типа второго рода исследования проводились лишь
при $n=0.$ Отметим работы [4,5].

В работе [6]  задача Коши для общего уравнения второго рода вида
$${{x}^{n}}{{V}_{xx}}-{{(-y)}^{m}}{{V}_{yy}}+a(x,y){{V}_{x}}+b(x,y){{V}_{y}}+c(x,y)V=f(x,y)\eqno(*)$$
решается сведением поставленной задачи к решению эквивалентного
интегро-дифференциального уравнения относительно искомого решения
$V(x,y)$.  В настоящем сообщении исследуется задача Коши хотя для
частного случая уравнения (*), т.е. при $a(x,y)=0,$ $b(x,y)=0,$
$c(x,y)=\mu ,$ $f(x,y)=0,$ но единственное решение задачи Коши для
уравнения (1) с условиями (3) находится методом Римана в явном
виде.

Предварительно установим ряд тождеств необходимых при решении
поставленной задачи.

Введем в рассмотрение ряды
$$F\left[
\begin{matrix}
   a,b  \\
   c;\,\,\sigma   \\
\end{matrix} \right]\equiv F(a,b;c;\sigma )=\sum\limits_{m=0}^{\infty }{\frac{{{(a)}_{m}}{{(b)}_{m}}}{m!{{(c)}_{m}}}\,}{{\sigma }^{m}},$$
$${{\Xi }_{2}}\left[ \begin{matrix}
   a,b  \\
   c;\,\,\sigma ,\rho   \\
\end{matrix} \right]\equiv {{\Xi }_{2}}(a,b;c;\sigma ,\rho )=\sum\limits_{m,k=0}^{\infty }{\frac{{{(a)}_{m}}{{(b)}_{m}}}{m!k!{{(c)}_{m+k}}}\,}{{\sigma }^{m}}{{\rho }^{k}},$$

$$\Phi \left[ \begin{matrix}
   a,b,c,d  \\
   e;\,\,\sigma ,\omega ,\rho   \\
\end{matrix} \right]\equiv {}_{3}\Phi _{B}^{(5)}\left[ \begin{matrix}
   a,b,c,d  \\
   e;\,\,\sigma ,\omega ,\rho   \\
\end{matrix} \right]=\sum\limits_{m,n,k=0}^{\infty }{\frac{{{(a)}_{m}}{{(b)}_{n}}{{(c)}_{m}}{{(d)}_{n}}}{m!n!k!{{(e)}_{m+n+k}}}\,}{{\sigma }^{m}}{{\omega }^{n}}{{\rho
}^{k}},\eqno(7)$$

$${{\Xi }_{pq}}\left[ \begin{matrix}
   a,b;\,{a}',{b}'  \\
   c;\,\,{c}',{d}';\sigma ,\rho   \\
\end{matrix} \right]=\sum\limits_{m,k=0}^{\infty }{\frac{{{(a)}_{m}}{{(b)}_{m}}{{({a}')}_{pk}}{{({b}')}_{q(m+k)}}}{m!k!{{(c)}_{m+k}}{{({c}')}_{pk}}{{({d}')}_{q(m+k)}}}\,}{{\sigma }^{m}}{{\rho
}^{k}},\eqno(8)$$ $${{\Psi }_{pq}}\left[ \begin{matrix}
   a,b,c,d;\,{a}',{b}'  \\
   c;\,\,{c}',{d}';\sigma ,\theta ,\rho   \\
\end{matrix} \right]=$$
$$=\sum\limits_{m,k=0}^{\infty
}{\frac{{{(a)}_{m}}{{(c)}_{m}}{{({a}')}_{pk}}{{({b}')}_{q(m+k)}}}{m!k!{{(e)}_{m+k}}{{({c}')}_{pk}}{{({d}')}_{q(m+k)}}}F\left[
\begin{matrix}
   b,e-d+m+k  \\
   e+m+k;\,\,\theta   \\
\end{matrix} \right]\,}{{\sigma }^{m}}{{\rho }^{k}},\eqno(9)$$
где  $a,b,c,d,e,{a}',{b}',{c}',{d}'$- комплексные параметры,
причем нижние параметры отличны от $0,\,-1,\,-2,...;$${{(\alpha
)}_{0}}=1,$ ${{(\alpha )}_{l}}=\alpha (\alpha +1)...(\alpha
+l-1),$   $l=1,2,....$ Здесь и далее $p$ и $q$ принимают значения
0 и 1; $\left| \sigma \right|<1,$ $\left| \omega
\right|<1,$$\theta =\omega /(\omega -1);$$F,\,$${{\Xi }_{2}}$ и
$\Phi $ - известные гипергеометрические функции [7,8].

$1^0.$ Справедливы соотношения $$\Xi_2\left[
\begin{matrix}
   a,b  \\
   c;\,\,\sigma ,\rho   \\
\end{matrix} \right]-{{\Xi }_{10}}\left[ \begin{matrix}
   a,b;\,e  \\
   c;\,\,e+1;\sigma ,\rho   \\
\end{matrix} \right]=$$
                                       $$=\frac{\rho }{(e+1)c}{{\Xi }_{10}}\left[ \begin{matrix}
   a,b;\,e+1  \\
   c+1;\,e+2;\sigma ,\rho   \\
\end{matrix} \right],\eqno(10)$$
$${{\Xi }_{10}}\left[ \begin{matrix}
   a,b;\,{a}'  \\
   c;\,\,{a}'+1;\sigma ,\rho   \\
\end{matrix} \right]-\frac{ab}{c(c+1)}\,\sigma \,{{\Xi }_{10}}\left[ \begin{matrix}
   a+1,b+1;\,{a}'  \\
   c+2;\,\,{a}'+1;\sigma ,\rho   \\
\end{matrix} \right]-$$
                $$-\frac{{a}'-c}{c(c+1)({a}'+1)}\rho \,{{\Xi }_{10}}\left[ \begin{matrix}
   a,b;\,{a}'+1  \\
   c+2;\,\,{a}'+2;\sigma ,\rho   \\
\end{matrix} \right]={{\Xi }_{2}}\left[ \begin{matrix}
   a,b  \\
   c+1;\,\sigma ,\rho   \\
\end{matrix} \right],\eqno(11)$$
$$(e-b-1){{\Psi }_{p1}}\left[ \begin{matrix}
   a,b,c,d;\,{a}',e-b  \\
   e;\,\,{c}',e-b-1;\sigma ,\theta ,\rho   \\
\end{matrix} \right]=(e-1){{\Psi }_{p0}}\left[ \begin{matrix}
   a,b,c,d-1;\,{a}'  \\
   e-1;\,\,{c}';\sigma ,\theta ,\rho   \\
\end{matrix} \right]-$$
                                         $$-b\,{{\Psi }_{p0}}\left[ \begin{matrix}
   a,b+1,c,d;\,{a}'  \\
   e;\,\,{c}';\sigma ,\theta ,\rho   \\
\end{matrix} \right],\eqno(12)$$
$$2(e-b-1)\,{{\Psi }_{p1}}\left[ \begin{matrix}
   a,b,c,d;\,{a}',e-b  \\
   e;\,\,{c}',e-b-1;\sigma ,\theta ,\rho   \\
\end{matrix} \right]+b\,(1-\theta )\,{{\Psi }_{p0}}\left[ \begin{matrix}
   a,b+1,c,d;\,{a}'  \\
   e;\,\,{c}';\sigma ,\theta ,\rho   \\
\end{matrix} \right]=$$
$$=(2e-b-d-2){{\Psi }_{p0}}\left[ \begin{matrix}
   a,b,c,d;\,{a}'  \\
   e;\,\,{c}';\sigma ,\theta ,\rho   \\
\end{matrix} \right]+d{{\Psi }_{p0}}\left[ \begin{matrix}
   a,b,c,d+1;\,{a}'  \\
   e;\,\,{c}';\sigma ,\theta ,\rho   \\
\end{matrix} \right]+$$
   $$+\frac{2ac}{e}\,\,\sigma \,{{\Psi }_{p0}}\left[ \begin{matrix}
   a+1,b,c+1,d;\,{a}'  \\
   e+1;\,\,{c}';\sigma ,\theta ,\rho   \\
\end{matrix} \right]+\frac{2{{({a}')}_{p}}}{e\,{{({c}')}_{p}}}\rho \,{{\Psi }_{p0}}\left[ \begin{matrix}
   a,b,c,d;\,{a}'+1  \\
   e+1;\,\,{c}'+1;\sigma ,\theta ,\rho   \\
\end{matrix} \right],\eqno(13)$$
$$\frac{(b-e-1)(d+1)}{e(e+1)}\,\,\theta \,{{\Psi }_{p1}}\left[
\begin{matrix}
   a,b,c,d+2;{a}',e-b+2  \\
   e+2;{c}',e-b+1;\sigma ,\theta ,\rho   \\
\end{matrix} \right]+$$
$$+(1-\theta )\,{{\Psi }_{p0}}\left[ \begin{matrix}
   a,b,c,d;{a}'  \\
   e;{c}';\sigma ,\theta ,\rho   \\
\end{matrix} \right]=\frac{b-d-1}{e}\,\,\theta \,{{\Psi }_{p0}}\left[ \begin{matrix}
   a,b,c,d+1;{a}'  \\
   e+1;{c}';\sigma ,\theta ,\rho   \\
\end{matrix} \right]+$$
                        $$+(1-\theta )\,{{\Psi }_{p1}}\left[ \begin{matrix}
   a,b,c,d+1;{a}',e+1  \\
   e+1;{c}',e;\sigma ,\theta ,\rho   \\
\end{matrix} \right],\eqno(14)$$
$$(e-1)(1-\theta )\,{{\Psi }_{p0}}\left[ \begin{matrix}
   a,b,c,d-1;{a}'  \\
   e-1;{c}';\sigma ,\theta ,\rho   \\
\end{matrix} \right]-(d-1)(1-\theta )\,{{\Psi }_{p0}}\left[ \begin{matrix}
   a,b,c,d;{a}'  \\
   e;{c}';\sigma ,\theta ,\rho   \\
\end{matrix} \right]+$$
        $$+(d-b){{\Psi }_{p0}}\left[ \begin{matrix}
   a,b,c,d;{a}'  \\
   e;{c}';\sigma ,\theta ,\rho   \\
\end{matrix} \right]=(e-b){{\Psi }_{p1}}\left[ \begin{matrix}
   a,b-1,c,d;{a}',e-b+1  \\
   e;{c}',e-b;\sigma ,\theta ,\rho   \\
\end{matrix} \right],  \eqno(15) $$
$${{\Psi }_{11}}\left[ \begin{matrix}
   a,b,c,d;{a}',{b}'+1  \\
   e;{a}'+1,{b}';\sigma ,\theta ,\rho   \\
\end{matrix} \right]=\frac{ac}{e{b}'}\,\sigma \,{{\Psi }_{10}}\left[ \begin{matrix}
   a+1,b,c+1,d;{a}'  \\
   e+1;{a}'+1;\sigma ,\theta ,\rho   \\
\end{matrix} \right]+$$
$$+\frac{{a}'-{b}'}{e({a}'+1){b}'}\,\,\rho \,{{\Psi }_{10}}\left[
\begin{matrix}
   a,b,c,d;{a}'+1  \\
   e+1;{a}'+2;\sigma ,\theta ,\rho   \\
\end{matrix} \right]+{{\Psi }_{00}}\left[ \begin{matrix}
   a,b,c,d  \\
   e;\sigma ,\theta ,\rho   \\
\end{matrix} \right].\eqno(16)$$

Доказательство равенств  (10) и (11) проводится путем разложения в
бесконечные ряды и сравнения коэффициентов.

Доказательство равенств (12)-(16) строится аналогично, исходя из
определения (9)  и известных соотношений для гипергеометрических
функций Гаусса [7, с.111].

$2^0.$  Известны следующие формулы дифференцирования
$$\frac{\partial }{\partial \sigma }\Phi =\frac{ac}{e}\Phi \left[
\begin{matrix}
   a+1,b,c+1,d  \\
   e+1;\sigma ,\omega ,\rho   \\
\end{matrix} \right],\eqno(17)$$
$$\frac{\partial }{\partial \omega }\Phi =\frac{bd}{e}\Phi \left[
\begin{matrix}
   a,b+1,c,d+1  \\
   e+1;\sigma ,\omega ,\rho   \\
\end{matrix} \right], \,\,\,\, \frac{\partial }{\partial \rho }\Phi
=\frac{1}{e}\Phi \left[ \begin{matrix}
   a,b,c,d  \\
   e+1;\sigma ,\omega ,\rho   \\
\end{matrix} \right].\eqno(18)$$

Используя формулу дифференцирования гипергеометрических функций
Гаусса [7], имеем
          $$\frac{\partial }{\partial \theta }\left[ {{(1-\theta )}^{b}}{{\Psi }_{pq}} \right]=-\frac{bd}{e}{{(1-\theta )}^{b-1}}{{\Psi }_{pq}}\left[ \begin{matrix}
   a,b+1,c,d+1;{a}',{b}'  \\
   e+1;{c}',{d}';\sigma ,\theta ,\rho   \\
\end{matrix} \right].\eqno(19)$$

$3^0.$  Из (9) при $\theta \to 1$ и $\operatorname{Re}(d-b)>0$
получаем формулу $${{\Psi }_{pq}}\left[ \begin{matrix}
   a,b,c,d;{a}',{b}'  \\
   e;{c}',{d}';\sigma ,1,\rho   \\
\end{matrix} \right]=\frac{\Gamma (e)\Gamma (d-b)}{\Gamma (e-b)\Gamma (d)}\,\,\,{{\Xi }_{pq}}\left[ \begin{matrix}
   a,c;{a}',{b}'  \\
   e-b;{c}',{d}';\sigma ,\rho   \\
\end{matrix} \right]. \eqno(20)$$

$4^0.$  Принимая во внимание формулу аналитического продолжения
гипергеометрических функций [7] нетрудно заключить, что функции
(7) и (9) связаны соотношением $$\Phi \left[
\begin{matrix}
   a,b,c,d  \\
   e;\sigma ,\omega ,\rho   \\
\end{matrix} \right]=\,\,{{(1-\omega )}^{-b}}\,{{\Psi }_{00}}\left[ \begin{matrix}
   a,b,c,d  \\
   e;\sigma ,\omega /(\omega -1),\rho   \\
\end{matrix} \right].\eqno(21)$$

Отсюда и из (20) простым вычислением получаем
($\operatorname{Re}(d-b)>0$) $$\underset{\omega \to
1}{\mathop{\lim }}\,\left\{ {{(1-\omega )}^{-b}}\,\Phi \left[
\begin{matrix}
   a,b,c,d  \\
   e;\sigma ,\omega /(\omega -1),\rho   \\
\end{matrix} \right] \right\}=\frac{\Gamma (e)\Gamma (d-b)}{\Gamma (e-b)\Gamma (d)}\,\,\,{{\Xi }_{2}}\left[ \begin{matrix}
   a,c  \\
   e;\sigma ,\rho   \\
\end{matrix} \right].\eqno(22)$$

$5^0.$ Известная формула автотрансформации [7]  дает формулу:
$${{\Psi }_{pq}}\left[ \begin{matrix}
   a,b,c,d;{a}',{b}'  \\
   e;{c}',{d}';\sigma ,\theta ,\rho   \\
\end{matrix} \right]={{(1-\theta )}^{d-b}}{{\Psi }_{pq}}\left[ \begin{matrix}
   a,d,c,b;{a}',{b}'  \\
   e;{c}',{d}';\sigma ,\theta ,\rho   \\
\end{matrix} \right].\eqno(23)$$

Переходим к решению задачи Коши.

Функция Римана для уравнения (1) известна [9]: $$R(\xi ,\eta
;{{\xi }_{0}},{{\eta }_{0}})={{\left( \frac{\eta +\xi }{{{\eta
}_{0}}+{{\xi }_{0}}} \right)}^{\alpha }}{{\left( \frac{\eta -\xi
}{{{\eta }_{0}}-{{\xi }_{0}}} \right)}^{\beta }}\Phi \left[
\begin{matrix}
   \alpha ,\beta ,1-\alpha ,1-\beta   \\
   1;\,\,\,\sigma ,\omega ,\rho   \\
\end{matrix} \right],\eqno(24)$$
где $$\sigma =\frac{({{\eta }_{0}}-\eta )(\xi -{{\xi
}_{0}})}{(\eta +\xi )({{\eta }_{0}}+{{\xi }_{0}})},   \omega
=\frac{(\eta -{{\eta }_{0}})(\xi -{{\xi }_{0}})}{(\eta -\xi
)({{\eta }_{0}}-{{\xi }_{0}})}, \rho =-{{\lambda }^{2}}({{\eta
}_{0}}-\eta )(\xi -{{\xi }_{0}}).$$

С помощью формулы аналитического продолжения (21), функцию Римана
(24) можно переписать в удобном для дальнейших исследований виде
$$R={{\left( \frac{\eta +\xi }{{{\eta }_{0}}+{{\xi }_{0}}}
\right)}^{\alpha }}\frac{{{({{\eta }_{0}}-\xi )}^{-\beta }}{{(\eta
-{{\xi }_{0}})}^{-\beta }}}{{{(\eta -\xi )}^{-2\beta }}}{{\Psi
}_{00}}\left[ \begin{matrix}
   \alpha ,\beta ,1-\alpha ,1-\beta   \\
   1;\,\,\,\sigma ,\theta ,\rho   \\
\end{matrix} \right],\eqno(25)$$
где $\theta =\omega /(\omega -1).$

Применяя метод Римана к области ${{\Delta }_{\varepsilon }}$,
ограниченной отрезками прямых  $M{Q_{\varepsilon }}:\xi ={{\xi
}_{0}},$ $M{{P}_{\varepsilon }}:\eta ={\eta _0}$ и
${{P}_{\varepsilon }}{{Q}_{\varepsilon }}:\eta =\xi +\varepsilon
,$ $\varepsilon >0,$   получим
   $$u({{\xi }_{0}},{{\eta }_{0}})=\frac{1}{2}{{(uR)}_{{{Q}_{\varepsilon }}}}+\frac{1}{2}{{(uR)}_{{{P}_{\varepsilon }}}}+{{J}_{1}}+{{J}_{2}},\eqno
(26)$$ где $$J_1=-\frac{1}{2}\int\limits_{{{\xi }_{0}}}^{{{\eta
}_{0}}-\varepsilon }{{{\left[ R\left( {{u}_{\xi }}-{{u}_{\eta }}
\right) \right]}_{\eta =\xi +\varepsilon }}}d\xi ,$$
$$J_2=\frac{1}{2}\int\limits_{{{\xi }_{0}}}^{{{\eta
}_{0}}-\varepsilon }{{{\left[ \left( {R_{\xi }}-{R_{\eta
}}+\frac{4\beta }{\eta -\xi }R \right)u \right]}_{\eta =\xi
+\varepsilon }}}d\xi .$$

Из формулы (25) следует, что функция $R$ на линии $\eta =\xi $
обращается  в бесконечность порядка $-2\beta $$(0<-2\beta <1).$
Стало быть, в отличие от вырождающихся гиперболических уравнений
первого рода, слагаемые в (26), содержащие $u({{\eta
}_{0}}-\varepsilon ,{{\eta }_{0}})$ и $u({{\xi }_{0}},{{\xi
}_{0}}+\varepsilon )$, не исчезают при $\varepsilon \to 0$ (и даже
более: являются бесконечно большими). Кроме того, ${{J}_{2}}$
переходит в расходящийся интеграл. Поэтому в начале преобразуем
подынтегральное выражение в ${{J}_{2}}$.

На основании формул дифференцирования (17)-(19) нетрудно получить
$$R_\xi -R_\eta +\frac{4\beta }{\eta -\xi }R={{\left( \frac{\eta
+\xi }{{{\eta }_{0}}+{{\xi }_{0}}} \right)}^{\alpha }}{{\left(
\frac{\eta -\xi }{{{\eta }_{0}}-{{\xi }_{0}}} \right)}^{\beta
}}I,$$ где $$I=-\beta (1-\beta ){{(1-\theta )}^{\beta -1}}\left(
\frac{\partial \theta }{\partial \xi }-\frac{\partial \theta
}{\partial \eta } \right){{\Psi }_{00}}\left[ \begin{matrix}
   \alpha ,1+\beta ,1-\alpha ,2-\beta   \\
   2;\,\,\,\sigma ,\theta ,\rho   \\
\end{matrix} \right]+$$
$$+\frac{2\beta }{\eta -\xi }\Phi \left[ \begin{matrix}
   \alpha ,\beta ,1-\alpha ,1-\beta   \\
   1;\,\,\,\sigma ,\omega ,\rho   \\
\end{matrix} \right]+{{\lambda }^{2}}(\eta -\xi -{{\eta }_{0}}+{{\xi }_{0}})\,\Phi \left[ \begin{matrix}
   \alpha ,\beta ,1-\alpha ,1-\beta   \\
   2;\,\,\,\sigma ,\omega ,\rho   \\
\end{matrix} \right]+$$

$$+\alpha (1-\alpha )\,\frac{{{\eta }_{0}}-{{\xi }_{0}}-\eta +\xi
}{({{\eta }_{0}}+{{\xi }_{0}})(\eta +\xi )}\Phi \left[
\begin{matrix}
   1+\alpha ,\beta ,2-\alpha ,1-\beta   \\
   2;\,\,\,\sigma ,\omega ,\rho   \\
\end{matrix} \right].$$

Используя тождество $$\frac{\partial \theta }{\partial \xi
}-\frac{\partial \theta }{\partial \eta }=\frac{{{\eta
}_{0}}-{{\xi }_{0}}}{({{\eta }_{0}}-\xi )(\eta -{{\xi
}_{0}})}\left[ 2-\frac{({{\eta }_{0}}-{{\xi }_{0}})(\eta -\xi
)}{({{\eta }_{0}}-\xi )(\eta -{{\xi }_{0}})}-\frac{{{(\eta -\xi
)}^{2}}}{({{\eta }_{0}}-\xi )(\eta -{{\xi }_{0}})} \right]$$
выделим слагаемые, имеющие конечный предел: $${{R}_{\xi
}}-{{R}_{\eta }}+\frac{4\beta }{\eta -\xi }R=T(\xi ,\eta ;{{\xi
}_{0}},{{\eta }_{0}})({{I}_{1}}+{{I}_{2}}+{{I}_{3}}),$$ где
$$T(\xi ,\eta ;{{\xi }_{0}},{{\eta }_{0}})={{\left( \frac{\eta
+\xi }{{{\eta }_{0}}+{{\xi }_{0}}} \right)}^{\alpha
}}\frac{{{({{\eta }_{0}}-\xi )}^{\beta }}{{(\eta -{{\xi
}_{0}})}^{\beta }}}{{{({{\eta }_{0}}-{{\xi }_{0}})}^{2\beta }}},$$

$${{I}_{1}}={{\lambda }^{2}}{{({{\eta }_{0}}-{{\xi }_{0}})}^{\beta
}}{{(\eta -\xi )}^{1+\beta }}{{({{\eta }_{0}}-\xi )}^{-\beta
}}{{(\eta -{{\xi }_{0}})}^{-\beta }}\Phi \left[ \begin{matrix}
   \alpha ,\beta ,1-\alpha ,1-\beta   \\
   2;\,\,\,\sigma ,\omega ,\rho   \\
\end{matrix} \right]+$$
$$+\beta (1-\beta ){{({{\eta }_{0}}-\xi {}_{0})}^{-1}}{{(1-\theta
)}^{1+2\beta }}{{\Psi }_{00}}\left[ \begin{matrix}
   \alpha ,1+\beta ,1-\alpha ,2-\beta   \\
   2;\,\,\,\sigma ,\theta ,\rho   \\
\end{matrix} \right]-$$
$$-\alpha (1-\alpha )\frac{\eta -\xi }{({{\eta }_{0}}+{{\xi
}_{0}})(\eta +\xi )}{{(1-\theta )}^{\beta }}\Phi \left[
\begin{matrix}
   1+\alpha ,\beta ,2-\alpha ,1-\beta   \\
   2;\,\,\,\sigma ,\omega ,\rho   \\
\end{matrix} \right],\eqno(27)$$
$${{I}_{2}}=-{{\lambda }^{2}}({{\eta }_{0}}-{{\xi
}_{0}}){{(1-\theta )}^{\beta }}\Phi \left[ \begin{matrix}
   \alpha ,\beta ,1-\alpha ,1-\beta   \\
   2;\,\,\,\sigma ,\omega ,\rho   \\
\end{matrix} \right]-\frac{2\beta ({{\eta }_{0}}-{{\xi }_{0}})}{({{\eta }_{0}}-\xi )(\eta -{{\xi }_{0}})}{{I}_{20}}+$$
$$+\beta (1-\beta ){{(\eta -\xi )}^{-1}}{{(1-\theta )}^{2}}{{\Psi
}_{00}}\left[ \begin{matrix}
   \alpha ,2-\beta ,1-\alpha ,1+\beta   \\
   2;\,\,\,\sigma ,\theta ,\rho   \\
\end{matrix} \right],$$
$${{I}_{20}}=(1-\beta ){{\Psi }_{00}}\left[ \begin{matrix}
   \alpha ,2-\beta ,1-\alpha ,1+\beta   \\
   2;\,\,\,\sigma ,\theta ,\rho   \\
\end{matrix} \right]-{{\Psi }_{00}}\left[ \begin{matrix}
   \alpha ,1-\beta ,1-\alpha ,\beta   \\
   1;\,\,\,\sigma ,\theta ,\rho   \\
\end{matrix} \right],$$
$${{I}_{3}}=\alpha (1-\alpha )\,\,\frac{{{\eta }_{0}}-{{\xi
}_{0}}}{({{\eta }_{0}}+{{\xi }_{0}})(\eta +\xi )}{{(1-\theta
)}^{2\beta }}{{\Psi }_{00}}\left[ \begin{matrix}
   1+\alpha ,\beta ,2-\alpha ,1-\beta   \\
   2;\,\,\,\sigma ,\theta ,\rho   \\
\end{matrix} \right].$$

Применяя последовательно ко второму и третьему слагаемым в
${{I}_{2}}$ формулы (12) и (13), после приведения подобных, имеем:
$${{I}_{2}}={{I}_{21}}+{{I}_{22}}+{{I}_{23}},$$ где
$${{I}_{21}}=-{{\lambda }^{2}}({{\eta }_{0}}-{{\xi
}_{0}}){{(1-\theta )}^{\beta }}\Phi \left[ \begin{matrix}
   \alpha ,\beta ,1-\alpha ,1-\beta   \\
   2;\,\,\,\sigma ,\omega ,\rho   \\
\end{matrix} \right]-$$
$$-{{\lambda }^{2}}\beta ({{\eta }_{0}}-{{\xi }_{0}})\theta
\,{{\Psi }_{00}}\left[ \begin{matrix}
   \alpha ,1-\beta ,1-\alpha ,1+\beta   \\
   3;\,\,\,\sigma ,\theta ,\rho   \\
\end{matrix} \right],$$
$${{I}_{22}}=\beta (1+\beta ){{(\eta -\xi )}^{-1}}(1-\theta
)\,{{\Psi }_{00}}\left[ \begin{matrix}
   \alpha ,1-\beta ,1-\alpha ,2+\beta   \\
   2;\,\,\,\sigma ,\theta ,\rho   \\
\end{matrix} \right],$$
$${{I}_{23}}=\alpha \beta (1-\alpha ){{(\eta -\xi )}^{-1}}(1-\theta
)\,\sigma \,\,{{\Psi }_{00}}\left[ \begin{matrix}
   1+\alpha ,1-\beta ,2-\alpha ,1+\beta   \\
   3;\,\,\,\sigma ,\theta ,\rho   \\
\end{matrix} \right].$$

Формулы (14) и (23) позволяют преобразовать ${{I}_{21}}$ в
выражение, имеющее конечный предел: $${{I}_{21}}=-{{\lambda
}^{2}}({{\eta }_{0}}-{{\xi }_{0}}){{(1-\theta )}^{1+2\beta
}}{{\Psi }_{01}}\left[
\begin{matrix}
   \alpha ,1+\beta ,1-\alpha ,1-\beta ;3  \\
   3;\,\,\,2;\,\,\,\sigma ,\theta ,\rho   \\
\end{matrix} \right]-$$
$$-\frac{1}{6}{{\lambda }^{2}}(1+\beta )(2+\beta )({{\eta
}_{0}}-{{\xi }_{0}})\,\theta \,{{\Psi }_{01}}\left[ \begin{matrix}
   \alpha ,1-\beta ,1-\alpha ,2+\beta ;3+\beta   \\
   4;\,\,\,2+\beta ;\,\,\,\sigma ,\theta ,\rho   \\
\end{matrix} \right].\eqno(28)$$

Рассмотрим теперь выражения ${{I}_{23}}$, ${{I}_{3}}$ и их сумму:
${{I}_{233}}={{I}_{23}}+{{I}_{3}}.$ Применяя к ${{I}_{3}}$ формулу
автотрансформации (23) и к ${{I}_{233}}$ формулу (15), получим
$${{I}_{233}}=\frac{1}{2}\alpha (1-\alpha )\,\,\frac{{{\eta
}_{0}}-{{\xi }_{0}}}{({{\eta }_{0}}+{{\xi }_{0}})(\eta +\xi
)}\,{{I}_{2331}},\eqno(29)$$ где $${{I}_{2331}}=-\beta {{(1-\theta
)}^{1+2\beta }}\,\,{{\Psi }_{00}}\left[ \begin{matrix}
   1+\alpha ,1+\beta ,2-\alpha ,1-\beta   \\
   3;\,\,\,\sigma ,\theta ,\rho   \\
\end{matrix} \right]+$$
$$+(2+\beta )\,\,{{\Psi }_{01}}\left[ \begin{matrix}
   1+\alpha ,-\beta ,2-\alpha ,1+\beta ;3+\beta   \\
   3;\,\,\,2+\beta ;\,\,\,\sigma ,\theta ,\rho   \\
\end{matrix} \right].$$

Введем вспомогательные функции $$\psi (\xi ,\eta )=\frac{(\eta
+\xi )({{\eta }_{0}}-{{\xi }_{0}})[\eta ({{\eta }_{0}}-\xi )-\xi
(\eta -{{\xi }_{0}})]}{{{[\eta ({{\eta }_{0}}-\xi )+\xi (\eta
-{{\xi }_{0}})]}^{2}}},$$ $$\varphi (\xi ,\eta )=T(\xi ,\eta
;{{\xi }_{0}},{{\eta }_{0}}){{\Psi }_{11}}\left[ \begin{matrix}
   \alpha ,-\beta ,1-\alpha ,1+\beta ;\beta ,1+\beta   \\
   1;\,\,\,1+\beta ,\beta ;\sigma ,\theta ,\rho   \\
\end{matrix} \right].$$

Пользуясь выражением производной $\left( \frac{\partial }{\partial
\xi }+\frac{\partial }{\partial \eta } \right)\varphi (\xi ,\eta
)$ на линии $\eta =\xi $ и формулой (16), можем
${{I}_{22}}$представить следующим образом:
$${{I}_{22}}={{I}_{221}}+\frac{1}{2}{{T}^{-1}}(\xi ,\eta ;{{\xi
}_{0}},{{\eta }_{0}})\,\,\psi (\xi ,\eta )d\varphi (\xi ,\eta ),$$
где $${{I}_{221}}=\beta (1+\beta )\,\kappa (\xi ,\eta )\,{{\Psi
}_{00}}\left[ \begin{matrix}
   \alpha ,1-\beta ,1-\alpha ,2+\beta   \\
   2;\,\,\sigma ,\theta ,\rho   \\
\end{matrix} \right]-$$
$$-\frac{\alpha }{\eta +\xi }\,\,\psi (\xi ,\eta )\,{{\Psi
}_{11}}\left[ \begin{matrix}
   \alpha ,-\beta ,1-\alpha ,1+\beta ;\beta ,1+\beta   \\
   1;\,\,\,1+\beta ,\beta ;\sigma ,\theta ,\rho   \\
\end{matrix} \right]+$$
$$+\frac{1}{2}{{\lambda }^{2}}({{\eta }_{0}}+{{\xi }_{0}}-\eta -\xi
)\,\,\psi (\xi ,\eta )\,\,{{\Psi }_{11}}\left[ \begin{matrix}
   \alpha ,-\beta ,1-\alpha ,1+\beta ;1+\beta ,2+\beta   \\
   2;\,\,\,2+\beta ,1+\beta ;\sigma ,\theta ,\rho   \\
\end{matrix} \right]-$$
$$-(1+\beta )\,\theta \,{{\kappa }^{+}}(\xi ,\eta )\psi (\xi ,\eta
)\,{{\Psi }_{10}}\left[ \begin{matrix}
   1+\alpha ,1-\beta ,2-\alpha ,1+\beta ;\beta   \\
   2;\,\,1+\beta ;\sigma ,\theta ,\rho   \\
\end{matrix} \right]+$$
$$+\,\,\,{{\kappa }^{-}}(\xi ,\eta )\,\psi (\xi ,\eta )\,{{\Psi
}_{11}}\left[ \begin{matrix}
   1+\alpha ,-\beta ,2-\alpha ,1+\beta ;\beta ,2+\beta   \\
   2;\,\,\,1+\beta ,1+\beta ;\sigma ,\theta ,\rho   \\
\end{matrix} \right],\eqno(30)$$
$$\kappa (\xi ,\eta )=\frac{{{\eta }_{0}}-{{\xi }_{0}}}{({{\eta
}_{0}}-\xi )(\eta -{{\xi }_{0}})}\left( 1-\frac{1}{2}\,\,\psi (\xi
,\eta )\,\,\frac{{{\eta }_{0}}+{{\xi }_{0}}-\eta -\xi }{{{\eta
}_{0}}-{{\xi }_{0}}} \right),$$ $${{\kappa }^{\pm }}(\xi ,\eta
)=\pm \frac{\alpha (1-\alpha )}{2}\frac{{{\eta }_{0}}\pm {{\xi
}_{0}}-\eta \mp \xi }{({{\eta }_{0}}+{{\xi }_{0}})(\eta +\xi )}.$$

Легко заметить, что $$\underset{\varepsilon \to 0}{\mathop{\lim
}}\,\psi (\xi ,\xi +\varepsilon )=\frac{2({{\eta }_{0}}+{{\xi
}_{0}}-2\xi )}{{{\eta }_{0}}-{{\xi }_{0}}},
\underset{\varepsilon \to 0}{\mathop{\lim }}\,\kappa (\xi ,\xi
+\varepsilon )=\frac{4}{{{\eta }_{0}}-{{\xi }_{0}}}.$$
Следовательно, выражение ${{I}_{211}}$ при $\eta \to \xi $
ограничено.

Таким образом, мы получили $${{R}_{\xi }}-{{R}_{\eta
}}+\frac{4\beta }{\eta -\xi }R=T(\xi ,\eta ;{{\xi }_{0}},{{\eta
}_{0}})\,A(\xi ,\eta ;{{\xi }_{0}},{{\eta }_{0}})+\frac{1}{2}\psi
(\xi ,\eta )\,d\varphi (\xi ,\eta ),$$ где  $A(\xi ,\eta ;{{\xi
}_{0}},{{\eta }_{0}})$ определяется из формул  (27)-(30), т.е.
$$A(\xi ,\eta ;{{\xi }_{0}},{{\eta
}_{0}})={{I}_{1}}+{{I}_{21}}+{{I}_{221}}+{{I}_{233}},$$ и имеет
конечный предел  при $\eta \to \xi $.

Подставляя полученный результат  в (26), найдем, что $$u({{\xi
}_{0}},{{\eta }_{0}})=\frac{1}{2}{{B}_{\varepsilon }}(\xi ,{{\xi
}_{0}},{{\eta }_{0}})-\frac{1}{2}\int\limits_{{{\xi
}_{0}}}^{{{\eta }_{0}}-\varepsilon }{{{\left[ \left( {{u}_{\xi
}}-{{u}_{\eta }} \right)R \right]}_{\eta =\xi +\varepsilon }}}d\xi
+$$
$$+\frac{1}{2}\int\limits_{{{\xi }_{0}}}^{{{\eta }_{0}}-\varepsilon }{{{\left[ u\,A \right]}_{\eta =\xi +\varepsilon }}}d\xi -\frac{1}{4}\int\limits_{{{\xi }_{0}}}^{{{\eta }_{0}}-\varepsilon }{\varphi }\,d(u\,\psi
),\eqno(31)$$ где $${{B}_{\varepsilon }}(\xi ,{{\xi }_{0}},{{\eta
}_{0}})=~{{(uR)}_{{{Q}_{\varepsilon
}}}}+{{(uR)}_{{{P}_{\varepsilon }}}}+\frac{1}{2}{{(u\,\varphi
\,\psi )}_{{{P}_{\varepsilon }}}}-\frac{1}{2}{{(u\,\varphi \,\psi
)}_{{{Q}_{\varepsilon }}}}$$.

Нетрудно проверить, что $$\underset{\varepsilon \to
0}{\mathop{\lim }}\,{{B}_{\varepsilon }}(\xi ,{{\xi }_{0}},{{\eta
}_{0}})=0.\eqno(32)$$

Выполняя в (31) предельный переход при $\varepsilon \to 0$,
учитывая (32), (20), (22) и используя соотношения (10)-(11), после
несложных  вычислений получим формулу решения задачи Коши (4)-(6):
$$u({{\xi }_{0}},{{\eta }_{0}})=(\eta_0+\xi_0)^{-\alpha}(\eta_0-\xi_0)^{-2\beta-1} \int\limits_{{{\xi }_{0}}}^{{{\eta
}_{0}}}{{{H}_{1}}({{\xi }_{0}},{{\eta }_{0}};\xi ;\lambda )\,}\tau
(\xi )d\xi +$$
$$+(\eta_0+\xi_0)^{-\alpha}(\eta_0-\xi_0)^{-2\beta-1} \int\limits_{{{\xi
}_{0}}}^{{{\eta }_{0}}}{{{H}_{2}}({{\xi }_{0}},{{\eta }_{0}};\xi
;\lambda )\,}{\tau }'(\xi )d\xi +$$ $$+ (\eta_0+\xi_0)^{-\alpha}
\int\limits_{{{\xi }_{0}}}^{{{\eta }_{0}}}{{{H}_{3}}({{\xi
}_{0}},{{\eta }_{0}};\xi ;\lambda )\,}\nu (\xi )d\xi ,\eqno(33)$$
где $${{H}_{1}}({{\xi }_{0}},{{\eta }_{0}};\xi ;\lambda )={{\gamma
}_{1}}(\eta_0-\xi)^\beta (\xi-\xi_0)^\beta \xi^\alpha \,S({{\xi
}_{0}},{{\eta }_{0}};\xi ,\lambda ),$$
$${{H}_{2}}({{\xi }_{0}},{{\eta }_{0}};\xi ;\lambda )=-{{\gamma
}_{1}} (\eta_0+\xi_0-2\xi)(\eta_0-\xi)^\beta (\xi-\xi_0)^\beta
\xi^\alpha\,\,\,{{\Xi }_{10}}\left[
\begin{matrix}
   \alpha ,1-\alpha ;\,\beta   \\
   \beta ;\,\,1+\beta ;{{\sigma }_{0}},{{\rho }_{0}}  \\
\end{matrix} \right],$$
$${{H}_{3}}({{\xi }_{0}},{{\eta }_{0}};\xi ;\lambda )=-{{\gamma
}_{2}}{{({{\eta }_{0}}-\xi )}^{-\beta }}{{(\xi -{{\xi
}_{0}})}^{-\beta }}\xi^\alpha \,{{\Xi }_{2}}\left[ \begin{matrix}
   \alpha ,1-\alpha   \\
   1-\beta ;{{\sigma }_{0}},{{\rho }_{0}}  \\
\end{matrix} \right],$$

$$S({{\xi }_{0}},{{\eta }_{0}};\xi ,\lambda )=2(1+2\beta )\Xi
-\frac{\alpha }{\xi }({{\eta }_{0}}+{{\xi }_{0}}-2\xi )\Xi
-\frac{\partial \Xi }{\partial {{\sigma }_{0}}}\frac{\partial
{{\sigma }_{0}}}{\partial \xi }+4{{\rho }_{0}}\frac{\partial \Xi
}{\partial {{\rho }_{0}}},$$ $${{\gamma }_{1}}=\frac{\Gamma
(1+2\beta )}{2^{1-\alpha}{{\Gamma }^{2}}(1+\beta )}, {{\gamma
}_{2}}={{[2(1-2\beta )]}^{2\beta }}{{2}^{\alpha -1}}\frac{\Gamma
(1-2\beta )}{{{\Gamma }^{2}}(1-\beta )},  {{\sigma
}_{0}}=\frac{({{\eta }_{0}}-\xi )(\xi -{{\xi }_{0}})}{2\xi
\,({{\eta }_{0}}+{{\xi }_{0}})}$$,     $$\Xi \equiv \,{{\Xi
}_{10}}\left[ \begin{matrix}
   \alpha ,1-\alpha ;\beta   \\
   \beta ;1+\beta ;{{\sigma }_{0}},{{\rho }_{0}}  \\
\end{matrix} \right], {{\rho }_{0}}=-{{\lambda }^{2}}({{\eta }_{0}}-\xi )(\xi -{{\xi }_{0}}).$$

Возвращаясь теперь к переменным $x$ и $y$, решение задачи Коши для
уравнения (1) с начальными данными (3) находим  в виде
$$V(x,y)=\int\limits_0^1 L_1(x,y;\zeta ;\mu )\,\tau
_1(t)d\zeta +$$ $$\int\limits_0^1L_2(x,y;\zeta ;\mu )\,\tau
'_1(t)d\zeta +\int\limits_0^1L_3(x,y;\zeta ;\mu )\nu _1(t)d\zeta
,\eqno(34)$$ где $$t=x_0+(2\zeta -1)y_0,
x_0=\frac{2}{2-n}{x^{(2-n)/2}},y_0=\frac{2}{2-m}{(-y)^{(2-m)/2}},$$
$$L_i(x,y;\zeta ;\mu
)=H_i(x_0-y_0,x_0+y_0,t;-4\lambda ^2), i=1,2,3.$$

Справедлива следующая

Теорема. Если $\tau _1(x)\in C^3[0,1]$ и $\nu_1(x)\in C^2[0,1]$,
то функция $V(x,y),$    определенная формулой (34), является
дважды непрерывно дифференцируемым, притом единственным, решением
задачи Коши для уравнения (1) с начальными данными (3) в области
$D.$

Доказательство. Единственность решения поставленной задачи
вытекает из самого способа получения формулы (34). В
справедливости остальных утверждений теоремы можно убедиться
непосредственным вычислением.

Литература

1.  Сахабиева Г.А. О краевой задаче Гурса для уравнения
гиперболического типа в трехмерном пространстве//.
Дифференциальные уравнения. Труды пединститутов РСФСР. Рязань.
1975, вып.6. С.200-207.

2.  Салахитдинов М.С., Исломов Б. О некоторых краевых задачах со
смещением для уравнения
$-{{(-y)}^{m}}{{V}_{yy}}+{{x}^{n}}{{V}_{xx}}+\mu V=0$.//
Неклассические уравнения математической физики и теории ветвления.
Ташкент, "Фан", 1988. С.24-34.

3.  Уринов А.К., Каримов Ш.Т. Краевые задачи со смещением для
вырождающегося гиперболического уравнения. // Труды института
Математики и компьютерных технологий. Вып.IV.,  Ашгабад, Ылым,
1995. С.159-164.

4. Евдокимов Ф.Ф. Задача Коши для уравнения
${{u}_{xx}}-{{(-y)}^{m}}{{u}_{yy}}-{{\lambda }^{2}}u=0$.//
Дифференциальные уравнения. Труды пединститутов РСФСР. Рязань.
1978, вып.12, С.45-50.

5.  Салахитдинов М.С., Эргашев Т.Г. Интегральное представление
обобщенного решения задачи Коши в классе $R_{2k}^{\lambda }$ для
одного уравнения гиперболического типа второго рода.// Узбекский
математический журнал, 1995, 1, С.67-75.

6.  Макаров И.А. Задача Коши для уравнения с двумя линиями
вырождения второго рода. Математическая физика. Куйбышев, 1976,
С.3-7.

7.  Бейтмен Г., Эрдейи А. Высшие трансцендентные функции. Т.1.
М.:Наука, 1973. 296 с.

8.  Jain R.N. The confluent hypergeometric functions of three
variables. Proc. Nat. Acad.Sci., India. 1966, vol. 36, No 2,
p.395-408.

9.  Андреев А.А., Волкодавов В.Ф., Шевченко Г.Н. О функциях
Римана. Дифференциальные уравнения. Труды пединститутов РСФСР.
Рязань. 1974, вып.4, С.25-31.

\end{document}